\documentclass{article}
\usepackage{graphics}
\usepackage{amsmath,amssymb,amsthm}
\begin{document}

\title{Stable Model of $X_0(125)$}
\author{Ken McMurdy}
\maketitle

\begin{abstract}
In this paper we determine the components in the stable model of $X_0(125)$ over ${\mathbb{C}}_5$ by constructing in the language of \cite{C2} an explicit semi-stable covering.
We then offer empirical data regarding the placement of certain CM $j$-invariants in the supersingular disk of $X(1)$ over ${\mathbb{C}}_5$ which suggests a moduli-theoretic interpretation for the components of the stable model.
The paper then concludes with a conjecture regarding the stable model of $X_0(p^3)$ for $p>3$, which is as yet unknown.
\end{abstract}

\newcommand{\Z}{\mathbb{Z}}
\newcommand{\Q}{\mathbb{Q}}
\newcommand{\C}{\mathbb{C}}
\newcommand{\F}{\mathbb{F}}
\newcommand{\N}{\mathbb{N}}
\newcommand{\End}{\text{End}}
\newcommand{\Disc}{\text{Disc}}
\newcommand{\Aut}{\text{Aut}}

\theoremstyle{plain}
\newtheorem{claim}{Claim}[subsection]
\newtheorem{lemma}{Lemma}[subsection]
\newtheorem*{proposition}{Proposition}
\newtheorem{conjecture}{Conjecture}[subsection]
\newtheorem{guess}{Guess}[subsection]
\theoremstyle{definition}
\newtheorem*{definition}{Definition}
\newtheorem*{note}{Note}
\newtheorem{example}{Example}[subsection]
\newtheorem*{remark}{Remark}

\begin{note} This is an unofficial version of the paper.
The definitive version has been published in the LMS Journal of Computation and Mathematics, Volume 7, pp. 21--36.
We are very grateful to the LMS for their careful refereeing, editing, and publication of the official version.
\end{note}

\section{Introduction}
The purpose of this paper is to begin to advance the work of Deligne-Rapoport, Katz-Mazur, Edixhoven, and others regarding the stable model of $X_0(p^n)$.
In the simplest case of $X_0(p)$ ($p\neq 2,3$), the minimal resolution over $\Z_p$ is always semi-stable and was explicitly described in \cite[VI.6.16]{DR}.
Edixhoven later worked out the minimal resolution of the Katz-Mazur model for $X_0(p^n)$ over $\Z_p$ (special case of \cite[1.4]{E}), but this model is never semi-stable when $n\geq 2$.
However, in the $n=2$ case, Edixhoven did go on to work out a semi-stable model for $X_0(p^2)$ over a finite extension of ${\mathbb{Z}}_p^{\text{unr}}$ (\cite[2.1,2.5]{E}).
Unfortunately, as Edixhoven states in the introduction, his methods do not generalize sufficiently to calculate the stable model when $n>2$ because of ``wild ramification.''

Recently, more progress toward understanding the stable model of $X_0(p^n)$ has been made by Coleman using a moduli-theoretic approach.
In \cite{C1} Coleman showed that for $p>3$ the ordinary region of $X_0(p^n)$ has exactly $2n$ connected components (over $\C_p$), whose reductions can be described using Igusa curves.
He also was able to give a moduli-theoretic interpretation of the horizontal components of Edixhoven's model for $X_0(p^2)$.
In particular, the points on these components correspond to pairs $(E,C)$ such that $E/C[p]$ is in the language of \cite[3.3]{B} too-supersingular, i.e. that $E/C[p]$ has no canonical subgroup.

In the first half of this paper we begin with the model for $X_0(125)$ constructed in \cite[4]{M}, and apply techniques of rigid analysis to construct a semi-stable covering.
From this point on, all such statements will always mean over $R_p$, the ring of integers in $\C_p$.
Paraphrasing \cite[2]{C2}, this means that we will find a finite set of disjoint affinoids $\{A_i\}$ with good reduction, such that the complement in $X_0(125)$ of $\cup A_i$ is simply the disjoint union of annuli.
By \cite[2.1]{C2} this is equivalent to determining a semi-stable model for the curve.
The goal of the remaining sections, then, is to provide a moduli-theoretic interpretation of the components that leads to a conjecture regarding the stable model of $X_0(p^3)$.
To do this we first map the affinoids in the semi-stable covering down to $X(1)$ via an appropriate moduli-theoretic map (using formulas from \cite{M}).
Then we show by explicit calculations that the images contain the $j$-invariants of curves with a certain type of complex multiplication.
Along with similar data for $p=7$ and $p=13$, this leads to two conjectures which could be considered the main results of the paper.
Conjecture \ref{Conjecture:MainConjecture1} concerns the distribution of CM curves inside a given $ss$ disk of $X(1)$, and Conjecture \ref{Conjecture:MainConjecture2} generalizes the stable model description of $X_0(125)$ to a description of $X_0(p^3)$.
These conjectures will be at least partially proven in an upcoming joint work with Robert Coleman.
\section{Stable Model of $X_0(125)$}\label{Section:StableModel}

In this section we determine the stable model for $X_0(125)$ (genus $8$) by essentially constructing in the language of \cite{C2} a semi-stable covering for the curve.
The initial model will be the one determined in \cite[4]{M}, namely the following system of equations.
\begin{multline}\label{Eq:C125p}
f_{125}^+(x,y)=y^4 - x^5 + 5xy^3 + 15x^2y^2 + 25x^3y + 25x^4 + 5y^3\\
    + 5xy^2 - 25x^3 + 15y^2 + 25x^2 + 25y - 25x + 25 = 0
\end{multline}
\begin{equation}\label{Eq:C125}
xu^2-yu+5=0
\end{equation}
As modular functions, the $q$-expansions at infinity of $x$, $y$, and $u$ can be expressed in terms of the Dedekind eta-function as follows.
$$u(q)=\frac{\eta(q)}{\eta(q^{25})}\qquad x(q)=\frac{u(q^5)}{u(q)}\qquad y(q)=u(q^5)+\frac{5}{u(q)}$$
In \cite[4]{M} it is shown that $x$ and $y$ are functions on $X_0(125)$ which are fixed by the Atkin-Lehner involution, $w_{125}$.
Therefore Equation \eqref{Eq:C125p} actually describes the genus two quotient curve, 
$$X_0(125)^+=X_0(125)/w_{125}.$$
Then $u$ is a third function on $X_0(125)$ (actually a pullback of a function on $X_0(25)$) which generates, by Equation \eqref{Eq:C125}, the degree $2$ extension from $X_0(125)^+$ up to $X_0(125)$.

In later sections we will want to interpret the various components of our model in moduli-theoretic terms in the hopes of conjecturing what happens in general (for $X_0(p^3)$, $p\neq 5$).
For this section, however, we will stick to bare-bones rigid analysis.
First we show that $X_0(125)^+$ has good reduction by finding an explicit good reduction model.
Then we show that the $10$ ramification points in the degree two extension lie in two equidistant sets of $5$, all within a wide-open annulus.
The result of this information is that the stable model of $X_0(125)$ has $5$ components, four of genus $2$ and one of genus $0$.
In particular, two components with the same (stable) reduction as $X_0(125)^+$ will be switched under the Atkin-Lehner involution.
The other two genus $2$ affinoids are fixed by Atkin-Lehner and have trivial quotients.

\subsection{Good Reduction of $X_0(125)^+$}\label{Sec:GoodReduction}

At the end of \cite{M} there is a proof that $X_0(125)^+$ has good reduction which is based on a hyperelliptic model for the genus $2$ curve.
In order to understand the extension up to $X_0(125)$, however, this hyperelliptic model is not optimal.
For this reason we will now offer a proof using a different model, one which will be more appropriate for the full analysis of $X_0(125)$.
Essentially, we will show that with respect to the parameter $y$ the affinoid $v_5(y)=3/4$ is a genus $2$ affinoid (in the sense of \cite[1]{M}) with good reduction.
Then in the section that follows we will see that the $10$ ramification points in the degree two extension up to $X_0(125)$ are actually separated from this affinoid by a wide open annulus.
This is a key point in understanding the stable model of $X_0(125)$.

\begin{claim}\label{Claim:GoodReduction}
The curve $X_0(125)^+$ has good reduction. 
More specifically, the affinoid described by $v_5(y)=3/4$ is a genus $2$ affinoid with good reduction.
\end{claim}
\begin{proof}
We begin by choosing any $r\in\C_5$ satisfying $r^5+25r-25$ and making the change of variables $x_0=x-r$.
Plugging into $f^+(x,y)$ we then obtain a polynomial $g^+(x_0,y)$ whose coefficients are polynomials in $r$ and are given below in Table \ref{x0-y-model}.
The entry in the $x_0^i$ row and $y^j$ column is the coefficient of $x_0^iy^j$.
\begin{table}[ht]
\begin{center}
\begin{tabular}{  c | c | c | c | c | c |}
        &  $y^4$  &  $y^3$  &  $y^2$  &  $y^1$  &  $y^0$ \\ \hline
$x_0^5$ &         &         &         &         & $-1$ \\ \hline
$x_0^4$ &         &         &         &         & $-5r+25$ \\ \hline
$x_0^3$ &         &         &         & $25$   & $-10r^2+100r-25$  \\ \hline
$x_0^2$ &         &         & $15$   & $75r$  & $-10r^3+150r^2-75r+25$ \\ \hline
$x_0^1$ &         & $5$    & $30r+5$ & $75r^2$ & $-5r^4+100r^3-75r^2+50r-25$ \\ \hline
$x_0^0$ &  $1$   & $5r+5$ & $15r^2+5r+15$ & $25r^3+25$ & $25r^4-25r^3+25r^2$ \\ \hline
\end{tabular}
\end{center}\caption{Coefficients of $g^+(x_0,y)=f^+(x_0+r,y)$}\label{x0-y-model}
\end{table}

When $v_5(y)=3/4$, the Newton Polygon for $g^+$ considered as a polynomial in $x_0$ shows that $v_5(x_0)=1/2$.
Consequently, there are three terms with minimal valuation which on their own (ignoring other terms) would form the equation
\begin{equation}\label{Eq:Affinoid_x125p}
x_0^5+25x_0=15y^2.
\end{equation}
This motivates a second change of variables which will result in a good reduction model.
We first choose any $\alpha$ and $\beta$ in $\C_5$ with $v_5(\alpha)=1/2$ and $v_5(\beta)=3/4$.
Then we make the change of variables, $\alpha x_1=x_0$ and $\beta y_1=y$.
To be completely precise, then, the final choice of model for the curve $X_0(125)^+$ in terms of the parameters $x_1$ and $y_1$ is the equation
$$\frac{1}{15\beta^2} f_{125}^+(\alpha x_1+r,\beta y_1)=0.$$
This equation has integral coefficients and reduces modulo the maximal ideal of $R_5$ to the following equation over ${\bar{\mathbb{F}}}_5$.
\begin{equation}\label{x125p_red}
y_1^2=\frac{\alpha^5}{15\beta^2}x_1\left(x_1^4+\frac{25}{\alpha^4}\right)
\end{equation}
From Equation \ref{x125p_red} we see that as claimed the affinoid $v_5(y)=3/4$ (equivalently $v_5(y_1)=0$) is a genus $2$ affinoid with good reduction.
In rigid terms we have shown that the entire region $v_5(y)<3/4$ is simply one residue disk, while $v_5(y)>3/4$ describes five residue disks.
\end{proof}

\subsection{Rigid Distribution of Ramified Points}

We now turn our attention to understanding what happens to the curve $X_0(125)^+$ in the quadratic extension up to $X_0(125)$.
In light of the previous section, it is crucial that we understand where the ramification points lie in relation to the affinoid $v_5(y)=3/4$.
Equation \eqref{Eq:C125} tells us that the ramification points satisfy $y^2=20x$.
Substituting this into Equation \eqref{Eq:C125p}, we obtain a polynomial in $y$ satisfied by the $y$ coordinates of these $10$ ramification points.
The valuations of the coefficients of this degree $10$ polynomial, say $p_{\text{ram}}(y)=\sum_{i=0}^{10} a_i y^i$ are as follows.
\begin{center}
\begin{tabular}{ | c | c | c | c | c | c | c | c | c | c | c | c | } \hline
$i$        & 10 &     9    & 8 & 7 & 6 & 5 & 4 & 3 & 2 & 1 & 0 \\ \hline
$v_5(a_i)$ & 0  & $\infty$ & 3 & 4 & 4 & 5 & 5 & 6 & 6 & 7 & 7 \\ \hline
\end{tabular}\end{center}
The Newton polygon for $p_\text{ram}(y)$ then tells us that $v_5(y)=7/10$ at each root. 
As an immediate consequence all $10$ ramification points lie outside of the affinoid $v_5(y)=3/4$.
We will now show that the region described by $1/2<v_5(y)<3/4$ is in fact a wide open annulus (not parameterized by $y$ though), and that these $10$ points lie in two equidistant sets of $5$ within the circle $v_5(y)=7/10$.
The stable model of $X_0(125)$ will then follow from an argument similar to the one used in \cite[4.6]{M}.

\begin{claim}\label{Claim:s_Annulus}
The region of $X_0(125)^+$ described by $1/2<v_5(y)<3/4$ is a wide-open annulus.
\end{claim}
\begin{proof}
The main idea here is that (on this region) the terms $x_0^5$ and $15y^2$ of $g^+(x_0,y)$ have minimal valuation (by Newton Polygons).
In other words, the curve is well-approximated by the equation $x_0^5=15y^2$. 
This suggests that we should be able to parameterize the region using the annulus 
$$A_s=\{\ s\ |\ 1/5<v_5(s)<1/4\ \}$$ 
and a map close to $x_0=s^2$, $y=s^5/\sqrt{15}$ (for a fixed square root of $15$).
This can be made precise by applying a souped-up version of Hensel's Lemma (see \cite[2.3]{M}).
First set $x_0=s^2$ exactly, and then consider the polynomial $h(y)=s^{-10}g^+(s^2,s^5y/\sqrt{15})$.
Note that the coefficients of $h(y)$ are integral-valued functions on $A_s$.
Also it is straightforward to check that $v_5(h(1))>0$ and $v_5(h'(1))=0$ everywhere on $A_s$.
Therefore there is a unique integral-valued function on $A_s$ which is a root of $h(y)$ close to $y=1$.
To explicitly parameterize the region, then, we simply take $y$ to be $s^5/\sqrt{15}$ times this root.

\begin{equation}\label{Eq:AsParam}
x_0=s^2\qquad y=\frac{s^5}{\sqrt{15}}(1+\text{smaller terms on $A_s$})
\end{equation}

It is immediate that this defines a map from $A_s$ to $X_0(125)^+$, that it is an injection, and that the image is contained in the region $1/2<v_5(y)<3/4$.
The onto argument is a little more subtle, though.
Let $(x_0,y)$ be any point satisfying $1/2<v_5(y)<3/4$.
The Newton polygon of $g^+$ as a polynomial in $x_0$ shows that $2/5<v_5(x_0)<1/2$.
But then the Newton polygon of $g^+$ as a polynomial in $y$ shows that only two of the four points with this $x_0$ coordinate satisfy the condition $1/2<v_5(y)<3/4$.
On the other hand, there are two points with this $x_0$ coordinate in the image of $A_s$ simply by taking $s$ to be either square root of $x_0$.
Therefore the map must be onto by this simple counting argument.
\end{proof}

Previously we had shown that the $10$ ramified points in the degree $2$ extension up to $X_0(125)$ lay in the region $v_5(y)=7/10$.
Now we know from Equation \eqref{Eq:AsParam} that this region is in fact a circle parameterized by $v_5(s)=6/25$.
What we would like to do now is analyze the geometry of these $10$ points with respect to the parameter $s$.
One way to do this would be to explicitly work out enough terms in the power series satisfied by the $10$ $s$ coordinates to determine the relative distances.
It is possible (and much easier), however, to infer the desired information from the relative distances with respect to $x$ and $y$ coordinates.
This is the approach that we will take in proving the following claim.
\begin{claim}\label{Claim:TwoDiscs}
The $10$ ramification points of $X_0(125)^+$ lie in two equidistant sets of $5$, all within the circle $v_5(s)=6/25$, such that
$$\forall i\neq j\quad v_5(s_i-s_j)=6/25\quad\text{or}\quad v_5(s_i-s_j)=13/50.$$
\end{claim}
\begin{proof}
First move one root of $p_{\text{ram}}(y)$ to $0$ by a simple translation and then look at the Newton Polygon of the resulting polynomial.
This polynomial has $4$ roots of valuation $8/10$ and $5$ of valuation $7/10$.
From this we learn that for a fixed $y$ coordinate of a ramified point, say $y_0$, we have $v_5(y-y_0)=7/10$ at $5$ of the ramified points and $v_5(y-y_0)=8/10$ at the remaining four.
Actually, we know a lot more because $y\approx s^5/\sqrt{15}$ and we are working $5$-adically.
Therefore, when $v_5(s_1)=v_5(s_2)=6/25$ we have
$$v_5(s_1-s_2)>6/25\ \Leftrightarrow\ v_5(s_1^5-s_2^5)>6/5\ \Leftrightarrow\ v_5(y_1-y_2)>7/10.$$
It follows that the $10$ points at least break up into two sets of $5$ under the relation $v_5(s_i-s_j)>6/25$, but it does {\em not} follow that within each subset the points are equidistant.
To determine this fact we now look at the geometry of the $x$ coordinates.

Similar to what was done for the $y$ coordinates of the $10$ ramified points, it is straightforward to calculate the polynomial satisfied by the $x$ coordinates.
Then moving one root to $0$ we find that the difference of any two $x$ coordinates has valuation precisely $1/2$.
Of course the same can then be said about the $x_0$ coordinates, since $x$ and $x_0$ differ by a constant.
So now choose any two ramified points, say with $s$ coordinates $s_1$ and $s_2$, satisfying $v_5(s_1-s_2)>6/25$.
We can determine exactly how close these points are from the following formula.
$$\tfrac{1}{2}=v_5(x_0(s_1)-x_0(s_2))=v_5(s_1^2-s_2^2)=v_5(s_1-s_2)+v_5(s_1+s_2)$$
The point is that 
$$v_5(s_1-s_2)>\tfrac{6}{25}\ \Rightarrow\ v_5(s_1+s_2)=\tfrac{6}{25}\ \Rightarrow\ v_5(s_1-s_2)=\tfrac{13}{50}.$$
This proves the remaining part of the claim.
\end{proof}
\subsection{Stable Model}
In this section we finally give a complete description of a semi-stable model of $X_0(125)$.
Essentially this will be done in two steps. 
First we show that lying over the two minimal affinoid disks of Claim \ref{Claim:TwoDiscs} (containing the ramification points) there are two genus $2$ affinoids in $X_0(125)$ with good reduction.
Since $X_0(125)^+$ has good reduction at $v_5(y)=3/4$, and the extension is a trivial two sheeted cover over this region, there will then be two more genus $2$ affinoids in $X_0(125)$ which are isomorphic copies of this one.
So just from the fact that the genus of $X_0(125)$ is $8$, there can be no other nontrivial components in the stable model.
However, in the main part of Claim \ref{Claim:StableModel} we show that there is one component of genus $0$ which meets the four genus $2$ affinoids in distinct residue classes.
This determines the stable model, and a picture summarizing this data along with the results of \cite[4]{M} then concludes the section.

\begin{claim}\label{Claim:GoodReduction2}
Let $s_1$ be the $s$ coordinate of any of the $10$ ramification points in the extension from $X_0(125)^+$ up to $X_0(125)$.
The region lying over the affinoid disk $v_5(s-s_1)\geq 13/50$ is a genus $2$ affinoid with good reduction.
\end{claim}
\begin{proof}
For convenience we reparameterize this disk with the disk described by $v_5(t)\geq 0$, by letting $t=\beta(s-s_1)$ where $\beta\in \C_5$ is anything satisfying $v_5(\beta)=-13/50$.
This places one ramification point at $t=0$ and the other $4$ on the circle $v_5(t)=0$.
To fix notation, say the $t$ coordinates of the other $4$ ramification points are $\alpha_1$, $\alpha_2$, $\alpha_3$, and $\alpha_4$, where $v_5(\alpha_i)=0$ and $v_5(\alpha_i-\alpha_j)=0$ for $i\neq j$.

Now, we know that the quadratic extension can be obtained by taking the square root of an appropriate analytic function.
For example, if we let $z=2xu-y$, Equation \eqref{Eq:C125} simply becomes $z^2=y^2-20x$.
The function $y^2-20x$ has been shown to have simple roots at the $5$ ramification points and no other roots (on this disk).
Therefore we may rewrite Equation \eqref{Eq:C125} in terms of $t$ and $z$ as
$$z^2=t(t-\alpha_1)(t-\alpha_2)(t-\alpha_3)(t-\alpha_4)P(t)$$
where $P(t)=a_0+a_1t+a_2t^2+\cdots$ is an analytic and nonvanishing function on the disk $v_5(t)\geq0$.
But this means that $v_5(a_0)<v_5(a_i)$ for all $i>0$ (see for example \cite[6.2.2]{R}).
Therefore by making the subsitution $z_1=z/\sqrt{a_0}$ we arrive at a model for the quadratic extension of the disk which has the following reduction.
\begin{equation}
z_1^2=t(t-\alpha_1)(t-\alpha_2)(t-\alpha_3)(t-\alpha_4)
\end{equation}
Since the roots of the right hand side are distinct this is the equation for a genus $2$ curve over $\bar{\F}_5$, which proves the claim.
\end{proof}

\begin{claim}\label{Claim:StableModel}
The stable model of $X_0(125)$ has $5$ components: $4$ genus $2$ components which do not intersect each other, and $1$ genus $0$ component which intersects each of the others in exactly one place.
\end{claim}
\begin{proof}
By using the quadratic formula and an appropriate expansion for the square root, it is straightforward to show that the regions $v_5(y)<7/10$ and $v_5(y)>7/10$ of $X_0(125)^+$ split trivially in the degree two extension (as in \cite[4.3]{M}).
Therefore over the affinoid $v_5(y)\geq3/4$, which was shown in \ref{Sec:GoodReduction} to have good reduction and genus $2$, we have two isomorphic copies in $X_0(125)$.
In addition, we have just shown in Claim \ref{Claim:GoodReduction2} that there are two more genus $2$ affinoids with good reduction lying over two disks within the circle $v_5(s)=6/25$.
Since the genus of the whole curve $X_0(125)$ is only $8$, there can be no other nontrivial components in the stable model.
All that remains, then, is to understand how these components fit together and whether or not there are also genus $0$ components in the stable model.

To answer this question, it suffices to look at the reduction of Equation \eqref{Eq:C125} over an appropriate affinoid.
So choose ramification points $s_1$ and $s_2$ as before with $v_5(s_1-s_2)=6/25$, and define an affinoid $B\subseteq A_s$ by
$$B=\{\ s\in A_s\ |\ v_5((s-s_1)(s-s_2))=12/25\ \}.$$
In other words, $B$ is just the circle $v_5(s)=6/25$ minus the two residue disks containing the ramification points in the extension up to $X_0(125)$.
Let $\hat{B}$ denote the affinoid of $X_0(125)$ lying over $B$.

Recall that on this region $v_5(x)=2/5$ and $v_5(x-r)>2/5$ for a particular $r\in\C_5$ which satisfies $r^5+25r-25=0$.
It follows then that if we make the substitutions $s_0=s/\alpha$ and $u_0=u/\beta$ where $v_5(\alpha)=6/25$ and $v_5(\beta)=3/10$, Equation \eqref{Eq:C125} reduces to the following.
\begin{equation}\label{Eq:Genus0Component}
u_0^2-\frac{\alpha^5}{\sqrt{15}\beta r}s_0^5 u_0 +\frac{5}{\beta^2 r}=0
\end{equation}
Equation \eqref{Eq:Genus0Component} is nonsingular over the reduction of $\hat{B}$.
Indeed, the two residue disks which were removed from $A_s$ correspond precisely to the $s_0$ coordinates of the two (finite) singular points over $\bar{\F}_5$.
Therefore $\hat{B}$ has good reduction.
Furthermore, the four genus $2$ affinoids meet the reduction of $\hat{B}$ in four distinct residue classes, the two singular points and the two distinct $s_0=0$ points.
Therefore this genus $0$ component can not be blown down and the stable model can only be as claimed.
\end{proof}
\begin{remark} It would be interesting to determine a precise field extension over which the stable model is defined, and the resulting action of the Galois group on the special fiber.
Unfortunately, this does {\em not} follow from our calculations.
Specifically, we do not have explicit equations for the two genus $2$ components determined in Claim \ref{Claim:GoodReduction2}.
What we {\em do} know is that the ramification index of any such field must be divisible by $100$, since the width of each of the four annuli bounding the genus $2$ components is in fact $1/100$.
\end{remark}
\begin{figure}[ht]
\begin{center}
\includegraphics{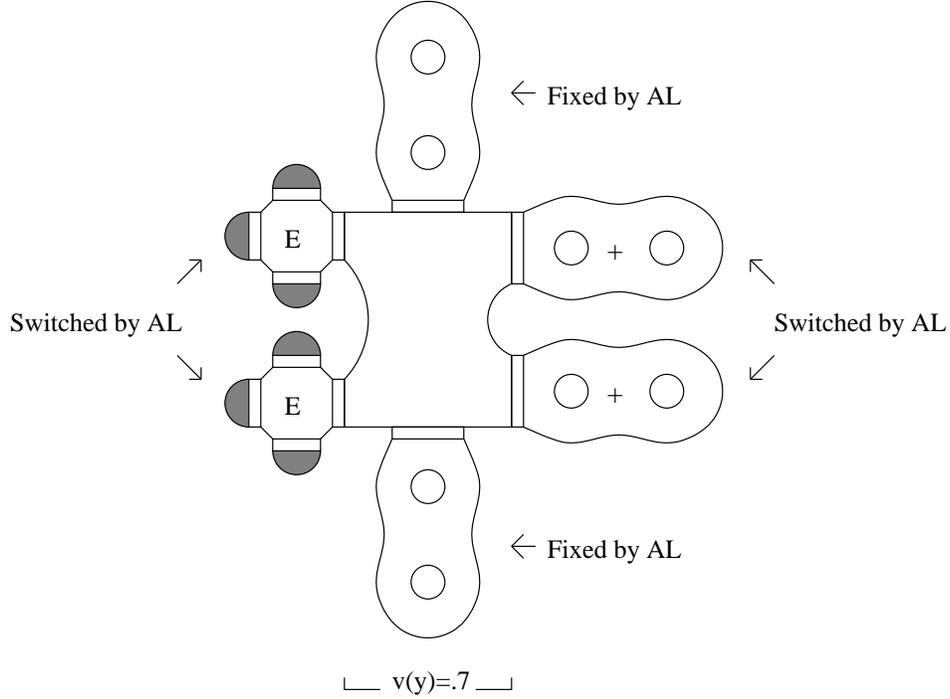}
\end{center}\caption{$X_0(125)/\C_5$, Semi-stable Covering}\label{Figure:Stable125}
\end{figure}
\begin{note} Figure \ref{Figure:Stable125} reflects all of the information from the preceeding claims regarding the stable model of $X_0(125)$. 
However, while the entire region described by $v_5(y)<7/10$ consists of two residue disks, these two residue disks were shown in \cite[4]{M} to contain components of great moduli-theoretic import.
In particular, this region contains six trivial ordinary components which are shaded in the picture.
We will also show in the following section that the two components marked ``E'' (the two components described by $v_5(y^2-5)=1$) map down via appropriate moduli-theoretic maps to the unique horizontal component of Edixhoven's model for $X_0(25)$.
Therefore, while not essential for a discussion of the stable model of $X_0(125)$, these components {\em are} essential for understanding the general conjecture regarding $X_0(p^3)$ which concludes the paper.
For this reason these features have also been included in the figure.
\end{note}

\section{Moduli-Theoretic Interpretation}
In this section we begin to formulate a moduli-theoretic interpretation of the components in the semi-stable model of $X_0(125)$ shown in Figure \ref{Figure:Stable125}.
This interpretation will then be the basis for a conjecture regarding the stable model of $X_0(p^3)$.
Philosophically, the main idea is to determine the image of each component in $X(1)$ via an appropriate map, and then ask what special moduli-theoretic properties are held by the elliptic curves corresponding to points in that region.
\subsection{Image in $X(1)$ of $X_0(125)$ Components}
To begin mapping the components of $X_0(125)$ down to $X(1)$, we first need to be very precise about which maps we are using.
\begin{definition} Let $M,N,d\in\N$ such that $dM|N$. Then we define a map $\pi_d$ from $X_0(N)$ to $X_0(M)$ in moduli-theoretic terms by
$$\pi_d(E,C)=(E/C[d],C[Md]/C[d]).$$
\end{definition}
\begin{note}
$\pi_d$ satisfies the compatibility condition, $\pi_{d_1}\circ \pi_{d_2}=\pi_{d_1d_2}$, wherever applicable.
For example, one could calculate $\pi_5:X_0(25)\to X(1)$ by way of either factorization, $\pi_5=\pi_5\circ\pi_1$ or $\pi_5=\pi_1\circ\pi_5$.
\end{note}
\noindent
Up until now we have simply used $u$ to refer to a parameter on $X_0(125)$, but in reality our $u$ is $\pi_1^*u$ where $u$ is a certain parameter on $X_0(25)$ and $\pi_1$ is the ``forgetful map'' as above.
So applying $\pi_1:X_0(125)\to X_0(25)$ amounts to just taking the $u$ coordinate.
We will not {\em need} a formula for $\pi_5:X_0(125)\to X_0(25)$.
For lower level moduli-theoretic maps, we will simply borrow the relevant formulas from \cite{M} and reproduce them for convenience in Table \ref{Table:Maps}.
\begin{table}
\begin{center}\renewcommand{\arraystretch}{1.5}
\begin{tabular}{| c | c |} \hline
 Moduli-Theoretic Map & Equation \\ \hline
      $\pi_1:X_0(5)\rightarrow X(1)$ & $\pi_1^*(j)=\frac{(t^2+2\cdot 5^3t+5^5)^3}{t^5}$ \\ 
      $\pi_5:X_0(5)\rightarrow X(1)$ & $\pi_5^*(j)=\frac{(t^2+10t+5)^3}{t}$ \\ 
        $w_5:X_0(5)\rightarrow X_0(5)$ & $w_5^*(t)=\frac{125}{t}$ \\ \hline
   $\pi_1:X_0(25)\rightarrow X_0(5)$ & $\pi_1^*(t)=\frac{u^5}{u^4+5u^3+15u^2+25u+25}$\\ 
   $\pi_5:X_0(25)\rightarrow X_0(5)$ & $\pi_5^*(t)=u(u^4+5u^3+15u^2+25u+25)$\\ 
 $w_{25}:X_0(25)\rightarrow X_0(25)$ & $w_{25}^*(u)=5/u$ \\ \hline
\end{tabular}
\end{center}\caption{Moduli-Theoretic Maps for Analysis of $X_0(5^n)$}\label{Table:Maps}
\end{table}

\begin{claim}\label{Claim:ImageE}
The two components marked ``E'' in Figure \ref{Figure:Stable125}, specifically the regions described by $v_5(y^2-5)=1$ are exactly the inverse images via $\pi_5$ and $\pi_{25}$ of the disk $v_5(j)\geq 5/2$ in $X(1)$.
\end{claim}
\begin{proof}
From the proof of \cite[4.4]{M}, it follows that one component, say $E_1$, is precisely the inverse image via $\pi_1$ of the affinoid in $X_0(25)$ described by $v_5(u^2-5)=1$.
The other component, say $E_2$, is the inverse image via $\pi_1$ of the circle $v_5(u)=1/10$.
To interpret $E_1$ we first apply $\pi_5:X_0(25)\to X_0(5)$ (using the formula from Table \ref{Table:Maps}) to see that $E_1$ is precisely the inverse image of the circle $v_5(t)=5/2$.
Then by composing with $\pi_1:X_0(5)\to X(1)$ we find that $E_1$ is actually the inverse image via $\pi_5$ of the disk $v_5(j)\geq 5/2$ inside $X(1)$.
Similarly we interpret $E_2$ by applying $\pi_{25}$.
First $E_2$ is the inverse image via $\pi_5$ of the circle $v_5(t)=1/2$ inside $X_0(5)$.
Then by composing with $\pi_5$ again we conclude that $E_2$ is actually the inverse image via $\pi_{25}$ of that same disk, $v_5(j)\geq 5/2$.
\end{proof}
\begin{claim}\label{Claim:ALImage} The region $v_5(y)\geq 7/10$ in $X_0(125)$ maps via $\pi_5$ onto the circle $v_5(j)=3/2$ in $X(1)$.
The images of the $4$ genus $2$ affinoids lie in $4$ distinct residue classes within this circle, specifically described by $v_5(j^4-5^6)>6$.
\end{claim}
\begin{proof}
From Equations \eqref{Eq:C125p} and \eqref{Eq:C125} it follows that $v_5(y)\geq 7/10$ is exactly the inverse image in $X_0(125)$ via $\pi_1$ of the region $v_5(u)=3/10$.
Then from the formula for $\pi_5:X_0(25)\to X_0(5)$, this is the entire inverse image of the circle $v_5(t)=3/2$ in $X_0(5)$.
By the formula for $\pi_1:X_0(5)\to X(1)$, this in turn maps onto the circle $v_5(j)=3/2$ inside $X(1)$.
Furthermore, since the maps are approximately $\pi_5^*t=u^5$ and $\pi_1^*j=t$ when restricted to the circles $v_5(u)=3/10$ and $v_5(t)=3/2$, the composition map is one-to-one on residue classes.
Now, in the proof of Claim \ref{Claim:StableModel} (see Equation \eqref{Eq:Genus0Component}) it was shown that $\pi_1$ maps the four genus $2$ components into four distinct residue disks within the circle $v_5(u)=3/10$.
Therefore, the only issue remaining is {\em which} four residue disks inside the circle $v_5(j)=3/2$ contain the images.

Recall that two of the affinoids contain the $10$ ramification points in the extension up from $X_0(125)^+$.
This set of points maps via $\pi_5$ onto the two points of $X_0(5)$ satisfying the equation $t^2=125$.
Therefore the images of these two affinoids in $X_0(5)$ lie in the two residue disks described by $v_5(t^2-125)>3$.
Of course this means that the affinoids map via $\pi_5$ into the residue disks $v_5(j^2-125)>3$ inside of $X(1)$.
Now we look at the other two affinoids, which by Equation \eqref{Eq:Genus0Component} land via $\pi_1$ inside the discs described by $v_5(u^2+5/r)>3/5$ where $r$ was a root of the equation $r^5+25r-25=0$.
Doing some quick arithmetic, this is equivalent to
$$v_5(u^{10}+5^5/r^5)>3\Leftrightarrow v_5(u^{10}+5^3)>3.$$
Then applying first $\pi_5:X_0(25)\to X_0(5)$ and then $\pi_1:X_0(5)\to X(1)$ these disks map onto the two disks of $X(1)$ described by $v_5(j^2+125)>3$.
Therefore, all four residue disks within the circle $v_5(j)=3/2$ can now be described by $v_5(j^4-5^6)>6$ as claimed.
\end{proof}
\begin{note}\label{Note:ALImage} From the formulas for $w_5$ and $\pi_1:X_0(5)\to X(1)$, the circle $v_5(t)=3/2$ is the Atkin-Lehner circle of $X_0(5)$ (fixed by $w_5$) and $v_5(j)=3/2$ is its image via the forgetful map.
\end{note}
\subsection{Computing the Too-Supersingular Region of $X(1)$}
In the previous section we defined two affinoids of $X_0(125)$, $E_1$ and $E_2$, which mapped via $\pi_5$ and $\pi_{25}$ respectively onto the disk $v_5(j)\geq 5/2$ inside of $X(1)$.
It is fairly straightforward to provide a moduli-theoretic description for the points lying in this disk.
In particular, we will now show that this disk consists of all points corresponding to an elliptic curve $E$ which is in the language of \cite{B} too-supersingular, i.e. a curve which has no canonical subgroup.
\begin{claim}\label{Claim:TooSS} $E/\C_5$ is too-supersingular iff $v_5(j(E))\geq 5/2$.
\end{claim}
\begin{proof}
Parameterize the supersingular disk using the disk $v_5(t)>0$ and the map which takes each $t$ to the $j$-invariant of the following curve.
$$E_t:y^2=x^3+tx+1\qquad j(t)=\frac{6912 t^3}{4t^3+27}$$
Strictly speaking this is a degree $3$ covering of the supersingular disk, ramified only at $j=0$.

Working out the (degree $12$) polynomial for the $x$ coordinates of $E_t[5]$, we see that the Newton Polygon has vertices $\{(0,0),(10,v_5(t)),(12,1)\}$ when $v_5(t)<5/6$, and  $\{(0,0),(12,1)\}$ otherwise.
If we take $z=x/y$ to be a parameter at infinity (containing all of $E_t[5]$ since $E_t$ is supersingular), this translates to the following information.
When $v_5(t)<5/6$ we have $v_5(z)=v_5(t)/20$ for $20$ of the points in $E_t[5]$, and 
$$v_5(z)=\tfrac{1-v_5(t)}{4}>\tfrac{v_5(t)}{20}$$ 
for the other four points (the canonical subgroup, along with $z=0$).
When $v_5(t)\geq 5/6$ we have $v_5(z)=1/24$ for all nonzero points of $E_t[5]$, and therefore there is no canonical subgroup.
This proves the claim since $v_5(j)=3v_5(t)$.
\end{proof}
\begin{note}
Robert Coleman has shown in \cite{C1} that the unique horizontal component (for a given $ss$ curve mod $p$) of Edixhoven's model for $X_0(p^2)$ is exactly the inverse image via $\pi_5$ of the corresponding too-ss disk.
This justifies the choice of names, $E_1$ and $E_2$, since it now follows from Claims \ref{Claim:ImageE} and \ref{Claim:TooSS} that these two components are simply $\pi_1^{-1}$ and $\pi_5^{-1}$ of Edixhoven's horizontal component for $X_0(25)$.
\end{note}
\subsection{Placement Data for CM Elliptic Curves in $X(1)$}
From Claims \ref{Claim:ImageE} and \ref{Claim:TooSS}, we have a clear way to interpret the components $E_1$ and $E_2$ in the semi-stable model for $X_0(125)$.
In particular, the component $E_1$ (resp. $E_2$) simply contains all points corresponding to pairs $(E,C)$ such that $E/C[5]$ (resp. $E/C[25]$) is too-ss.
Similarly, there is a clear moduli-theoretic interpretation for the union of the components of Claim \ref{Claim:ALImage}, as the entire region $v_5(y)\geq 7/10$ has been shown to be precisely $\pi_5^{-1}$ of the Atkin-Lehner circle of $X_0(5)$.
What is still unclear at this time is how to interpret the four special disks inside of the circle $v_5(j)=3/2$, described by $v_5(j^4-5^6)>6$, which contain the images of the four genus $2$ affinoids of $X_0(125)$.
In this section, we attempt to answer that question with empirical data placing certain types of ``CM points'' in those disks.
We begin by making a conjecture regarding the placement of the $j$-invariants of CM curves $E/\C_5$ such that $\End(E)\otimes \Z_5$ is the maximal order in a ramified (quadratic) extension of $\Q_5$.
\begin{conjecture}\label{Conjecture:CMp5} Suppose $E/\C_5$ is an elliptic curve with CM.\\
\indent (1) If $\End(E)\otimes \Z_5\cong \Z_5[\sqrt{-5}]$ then $v_5(j(E)^2-125)>3$.\\
\indent (2) If $\End(E)\otimes\Z_5\cong \Z_5[\sqrt{-10}]$ then $v_5(j(E)^2+125)>3$.
\end{conjecture}
To provide empirical evidence in support of the conjecture, we have worked out explicitly the $j$-invariants of various curves with these two types of CM.
For a given endomorphism ring $R$, this is done by first determining explicit representatives for the ideal class group of $R$ (using \cite[8 \S1]{L} for non-maximal orders).
Then writing each representative as $c(\Z+\tau\Z)$ we use the usual $q$-expansion formula to approximate $j(\tau)$ sufficiently well.
From the theory of CM curves these $j$-invariants (for a fixed $R$) are conjugate algebraic integers, so that the polynomial with these roots is monic and has integer coefficients.
The conjecture is verified then if this polynomial in $j$ is sufficiently close $5$-adically to a power of $j^2-125$ (Case 1) or $j^2+125$ (Case 2).
For easy verification of a few of the most basic examples (and to make the conjecture more concrete), we have included some complete sets of $\tau$ values in Tables \ref{Table:Case1Examples} and \ref{Table:Case2Examples}.
\begin{table}
\begin{center}\renewcommand{\arraystretch}{1.5}
\begin{tabular}{ | c | c | } \hline
$\End(E)$ & Values of $\tau$ \\ \hline
$\Z[\sqrt{-5}]$ & $\tau=\sqrt{-5}$, $\tfrac{1+\sqrt{-5}}{2}$ \\ \hline
$\Z[2\sqrt{-5}]$ & $\tau=2\sqrt{-5}$, $\tfrac{2(1+\sqrt{-5})}{3}$, $\tfrac{2(2+\sqrt{-5})}{3}$, $\tfrac{2\sqrt{-5}}{5}$ \\ \hline
$\Z[3\sqrt{-5}]$ &  $\tau=3\sqrt{-5}$, $\tfrac{3(1+\sqrt{-5})}{2}$, $\tfrac{3\sqrt{-5}}{5}$, $\tfrac{3(3+\sqrt{-5})}{7}$ \\ \hline
$\Z[\sqrt{-30}]$ & $\tau=\sqrt{-30}$, $\tfrac{\sqrt{-30}}{2}$, $\tfrac{\sqrt{-30}}{3}$, $\tfrac{\sqrt{-30}}{5}$ \\ \hline 
$\Z[\tfrac{1+\sqrt{-55}}{2}]$ & $\tau=\tfrac{1+\sqrt{-55}}{2}$, $\tfrac{1+\sqrt{-55}}{4}$, $\tfrac{-1+\sqrt{-55}}{4}$, $\tfrac{5+\sqrt{-55}}{10}$ \\ \hline   
$\Z[\sqrt{-70}]$ & $\tau=\sqrt{-70}$, $\tfrac{\sqrt{-70}}{2}$, $\tfrac{\sqrt{-70}}{5}$, $\tfrac{\sqrt{-70}}{7}$ \\ \hline   
\end{tabular}
\end{center}\caption{Examples for Case 1 of Conjecture \ref{Conjecture:CMp5}}\label{Table:Case1Examples}
\end{table}
\begin{table}
\begin{center}\renewcommand{\arraystretch}{1.5}
\begin{tabular}{ | c | c | } \hline
$\End(E)$ & Values of $\tau$ \\ \hline
$\Z[\sqrt{-10}]$ & $\tau=\sqrt{-10}$, $\tfrac{1}{2}\sqrt{-10}$ \\ \hline
$\Z[2\sqrt{-10}]$ & $\tau=2\sqrt{-10}$, $\tfrac{2\sqrt{-10}}{5}$, $\tfrac{2(2+\sqrt{-10})}{7}$, $\tfrac{2(1+\sqrt{-10})}{11}$ \\ \hline
$\Z[\tfrac{1+\sqrt{-15}}{2}]$ & $\tau=\tfrac{1+\sqrt{-15}}{2}$, $\tfrac{1+\sqrt{-15}}{4}$ \\ \hline
$\Z[\sqrt{-15}]$ & $\tau=\sqrt{-15}$, $\tfrac{\sqrt{-15}}{3}$\\ \hline
$\Z[\tfrac{1+\sqrt{-35}}{2}]$ & $\tau=\tfrac{1+\sqrt{-35}}{2}$, $\tfrac{5+\sqrt{-35}}{6}$\\ \hline
$\Z[\sqrt{-65}]$ & $\tau=\sqrt{-65}$, $\tfrac{1+\sqrt{-65}}{2}$, $\tfrac{1+\sqrt{-65}}{3}$, $\tfrac{-1+\sqrt{-65}}{3},$\\ 
                 & $\tfrac{\sqrt{-65}}{5}$, $\tfrac{1+\sqrt{-65}}{6}$, $\tfrac{-1+\sqrt{-65}}{6}$, $\tfrac{5+\sqrt{-65}}{10}$ \\ \hline
\end{tabular}
\end{center}\caption{Examples for Case 2 of Conjecture \ref{Conjecture:CMp5}}\label{Table:Case2Examples}
\end{table}

If true, the preceeding conjecture would provide at least the beginning of a moduli-theoretic description of the components in the stable model for $X_0(125)$.
However, at the end of the paper we will seek to propose a general conjecture for the stable model of $X_0(p^3)$.
Therefore, we will want to know how Conjecture \ref{Conjecture:CMp5} generalizes to other primes.
For this reason we have also checked a number of examples for $p=7$ and $p=13$, and have found that the data consistently supports Conjectures \ref{Conjecture:CMp7} and \ref{Conjecture:CMp13}. 
All three conjectures will be tied together with a general conjecture in next and final section of the paper.
\begin{conjecture}\label{Conjecture:CMp7} Suppose $E/\C_7$ is an elliptic curve with CM. Let 
$$j_0(E)=j(E)-1728.$$
\indent(1) If $End(E)\otimes\Z_7\cong \Z_7[\sqrt{-7}]$ then $v_7(j_0(E)^4-7^4)>4$.\\
\indent(2) If $End(E)\otimes\Z_7\cong \Z_7[\sqrt{-21}]$  then $v_7(j_0(E)^4+7^4)>4$.
\end{conjecture}
\begin{conjecture}\label{Conjecture:CMp13} Suppose $E/\C_{13}$ is an elliptic curve with CM. Let
$$j_0(E)=j(E)-5.$$
\indent(1) If $End(E)\otimes\Z_{13}\cong\Z_{13}[\sqrt{-13}]$ then $v_{13}(j_0(E)^{14}-13^7)>7$.\\
\indent(2) If $End(E)\otimes\Z_{13}\cong\Z_{13}[\sqrt{-26}]$ then $v_{13}(j_0(E)^{14}+13^7)>7$.
\end{conjecture}
\subsection{Conjectural Moduli-Theoretic Interpretation}
For a conjecture regarding the components in the stable model of $X_0(p^3)$, it is at least clear what ``types'' of components are suggested by the $X_0(125)$ data.
We will naturally conjecture that in general, for each supersingular elliptic curve (at least defined over $\F_p$), we have components that look like $E_1$, $E_2$, the four genus $2$ components, and the trivial component which intersects the others.
In other words, we will conjecture that there are components which have the same moduli-theoretic properties that the components in our semi-stable covering of $X_0(125)$ have been proven or conjectured to have.
Of the various ways in which to improve such a conjecture, one way would certainly be to conjecture the number of such components and give an explanation of the number.
This we will do with the help of Lemma \ref{Lemma:Ap}.
Conjecture \ref{Conjecture:MainConjecture1} uses the lemma to generalize the conjectures of the previous section regarding the placement of CM $j$-invariants.
Conjecture \ref{Conjecture:MainConjecture2} then takes this into account in describing the stable model of $X_0(p^3)$.
Finally, we conclude the paper with a guess which goes one step farther and predicts the genera of the components.
\begin{definition}
Choose $\alpha\in\F_p$ a quadratic non-residue, and $i\in \F_{p^2}$ with $i^2=\alpha$.
Then we define a finite $\F_p$-algebra
$$\bar{A}_p=\F_p[i,\epsilon_j,\epsilon_k]$$
where $\epsilon_j\epsilon_k=\epsilon_k\epsilon_j=\epsilon_j^2=\epsilon_k^2=0$ and $i\epsilon_j=\epsilon_k=-\epsilon_j i$.
\end{definition}
\begin{lemma}\label{Lemma:Ap} $\F_{p^2}^*$ acts on the nilradical ${\mathcal{N}}(\bar{A}_p)=\{c \epsilon_j+d\epsilon_k\}$ so that\\
\indent(1) $Sta(x)=\F_p^*$ and consequently $|Orb(x)|=p+1$ for all $x\neq 0$.\\
\indent(2) For all $x_1=c_1 \epsilon_j + d_1 \epsilon_k$ and $x_2=c_2 \epsilon_j + d_2 \epsilon_k$,
$$x_1\sim x_2\ \Leftrightarrow\ c_1^2-\alpha d_1^2=c_2^2-\alpha d_2^2$$
\end{lemma}
\begin{proof}The action is simply conjugation, after identifying $\F_{p^2}$ with $\F_p[i]$ of course.
Now, for part (1) let $x=c\epsilon_j+d\epsilon_k$ and $a+bi\in Sta(x)$.
Then by definition,
$$(a+bi)(c\epsilon_j+d\epsilon_k)\frac{a-bi}{a^2-\alpha b^2}=c\epsilon_j+d\epsilon_k.$$
Solving and setting $\epsilon_j$ and $\epsilon_k$ coefficients equal, this leads to the following system of linear equations in $c$ and $d$.
\begin{align*}
            (2b^2\alpha)c+(2ab\alpha)d&=0\\
                   (2ab)c+(2b^2\alpha)d&=0
\end{align*}
When $b=0$, i.e. when $a+bi\in \F_p^*$, the system is trivially satisfied for all $(c,d)$.
When $b\neq0$, however, the determinant is $4b^2\alpha(b^2\alpha-a^2)\neq0$ and therefore $x=0$ is the only solution.

Similarly, it is a straightforward exercise in congruences to verify half of Part (2), namely that $x_1\sim x_2$ implies $c_1^2-\alpha d_1^2=c_2^2-\alpha d_2^2$.
This means that as soon as $c^2-\alpha d^2=k\neq0$ has one solution, it must have at least $p+1$ solutions by part (1).
But $c^2-\alpha d^2=k$ always has a solution since the set $\{\alpha d^2+k\}$ has $(p+1)/2$ elements and therefore must contain a quadratic residue.
Therefore we are done since there are $p-1$ choices for $k\neq0$, $p+1$ solutions to $c^2-\alpha d^2=k$, and only $p^2-1$ nonzero elements of ${\mathcal{N}}(\bar{A}_p)-\{0\}$ to begin with.
\end{proof}
The reason that this lemma is relevant is the following.
Let $A$ be a supersingular curve so that $\End(A)\otimes\F_p\cong \bar{A}_p$.
The isomorphism is non-canonical, but it suffices to choose an isomorphism once and for all.
Now suppose that $E/\C_p$ is a CM curve with $p||\Disc(\End(E))$ and such that $\bar{E}\cong A$.
Identify $\End(E)$ with a subring of $\C_p$ via the canonical embedding coming from the action of $\End(E)$ on holomorphic differentials.
Then for every isomorphism $\lambda:\bar{E}\to A$ we obtain an embedding $\sigma_\lambda:\End(E)\to\End(A)$.
Furthermore, part (2) of Lemma \ref{Lemma:Ap} implies that for any uniformizer $u\in\End(E)$, the conjugacy class of $\sigma_\lambda(u)$ inside $\End(A)\otimes\F_p\cong \bar{A}_p$ is independent of $\lambda$.
So without taking $\lambda$ into account, there are only $p+1$ options for the image of $\sigma_\lambda(u)$ inside $\bar{A}_p$.
To illustrate this point, consider the following example.

\begin{example}\label{Example:p7} Let $E/\C_7$ be a curve with $\End(E)=\Z[\sqrt{-7}]$, and let $A$ be the unique supersingular elliptic curve in characteristic $7$.
From \cite[5.1]{P} we may take $\End(A)$ to be a maximal order in the quaternion algebra $\Q[i,j,k]$ with $i^2=-1$, $j^2=-7$, and $ij=-ji=k$.
Let $u$ be the uniformizer $2\sqrt{-7}\in \End(E)$ and suppose that for a particular isomorphism $\lambda:\bar{E}\to A$ we have $\sigma_\lambda(u)=a+bi+cj+dk$.
Since
$$(a+bi+cj+dk)^2=a^2-b^2-7c^2-7d^2+2a(bi+cj+dk)=-4\cdot7$$
we must immediately have $a=0$ and $b\equiv 0$ (mod $7$).
Therefore we have $c^2+d^2\equiv 4$, which means that for this example the image of $\sigma_\lambda(u)$ inside $\bar{A}_7$ must lie in the conjugacy class (of $8$ elements) containing
$$\{\pm 2\epsilon_j,\pm 2\epsilon_k,\pm3\epsilon_j\pm3\epsilon_k\}.$$
\end{example}

However, the precise image of $\sigma_\lambda(u)$ is {\em not} independent of $\lambda$.
So we can not yet use this image to put an equivalence relation on the set of curves with a given endomorphism ring.
To obtain an invariant which is independent of $\lambda$, we first note that $\Aut(A)$ acts transitively on the set of isomorphisms from $\End(\bar{E})$ to $\End(A)$ via conjugation, with $\pm1$ as the kernel in all cases.
So letting $i(A)=|\Aut(A)|/2$ we see that the conjugacy class defined by $u$ inside of $\bar{A}_p$ can be broken down into $(p+1)/i$ subsets, so that the {\em subset} containing the image of $\sigma_\lambda(u)$ is now independent of $\lambda$.
Again, we illustrate this point by revisiting the previous example.

\begin{example}\label{Example:p7b} The unique supersingular curve $A$ in characteristic $7$ has $|\Aut(A)|=4$.
In terms of the above identification, the automorphism group is generated by the (invertible) element $i$.
So for example, if $\sigma_\lambda(u)\equiv 2j$ for some $\lambda$, then $\sigma_{i\lambda}(u)\equiv i(2j)i^{-1}=-2j$.
This means that if we want to associate to $u$ an element of $\bar{A}_7$ which is {\em independent} of $\lambda$, we can {\em not} distinguish between $\pm 2\epsilon_j$.
Similarly, conjugation by $i$ breaks down the entire set of eight elements into four subsets of order $2$, namely
$$\{\pm 2\epsilon_j\},\{\pm 2\epsilon_k\},\{\pm(3\epsilon_j-3\epsilon_k)\},\ \text{and}\ \{\pm(3\epsilon_j+3\epsilon_k)\}.$$
So independent of $\lambda$, it makes sense to say that there are $4=(7+1)/2$ options for the image of $\sigma_\lambda(u)$ in $\bar{A}_7$.
\end{example}

The result of the preceeding argument is that we have shown how to put an equivalence relation with $(p+1)/i$ (possibly empty) classes on the set of CM curves $E$ reducing to a fixed supersingular curve $A$ and with a fixed endomorphism ring such that $p||\Disc(\End(E))$.
Indeed, with the same argument we may even generalize the relation by requiring only that $\End(E)\otimes\Z_p$ be fixed.
Since there are only two ramified quadratic extensions of $\Q_p$, this makes a total of $2(p+1)/i$ classes into which all such CM curves must fall.
What we have {\em not} shown, however, is that this equivalence relation is reflected somehow in the geometric placement of the $j$-invariants of these CM curves inside of $X(1)$.
So we have provided evidence for Conjecture \ref{Conjecture:MainConjecture1}, which would explain Conjectures \ref{Conjecture:CMp5}, \ref{Conjecture:CMp7}, and \ref{Conjecture:CMp13} with a general theory.
However, we do not yet have a proof of the result at this time.

\begin{figure}\begin{center}
\includegraphics{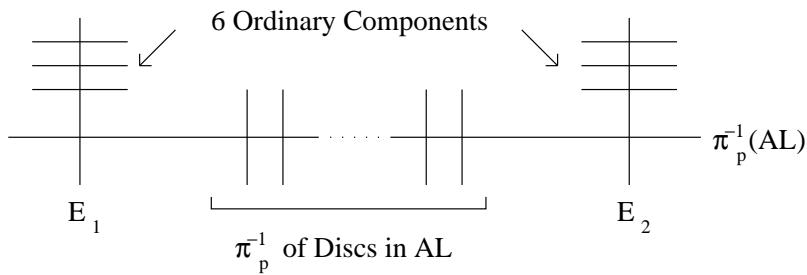}
\end{center}\caption{Conjectural Partial Graph of $X_0(p^3)$ Stable Reduction}\label{Figure:Stablep3Graph}
\end{figure}

\begin{conjecture}\label{Conjecture:MainConjecture1} Let $A$ be a supersingular elliptic curve defined over $\F_p$ and let $i=i(A)=|Aut(A)|/2$.\\
\indent (1) The $j$-invariants of all CM curves $E/\C_p$ such that $\bar{E}=A$ and 
$$p||\Disc(\End(E))$$
\indent\indent lie in $2(p+1)/i$ residue disks inside the corresponding Atkin-Lehner\\
\indent\indent circle, $(p+1)/i$ for each ramified quadratic extension of $\Q_p$.\\
\\
\indent (2) Two such curves, $E_1$ and $E_2$, lie in the same residue disk iff 
$$End(E_1)\otimes\Z_p=End(E_2)\otimes\Z_p$$
\indent\indent and for any uniformizer $u\in\C_p$ of the common image we have\\
$$\sigma_{\lambda_1}(u)\equiv\sigma_{\lambda_2}(u)\in\bar{A}_p$$
\indent\indent for some $\lambda_i:\bar{E}_i\to A$, $i=1,2$.
\end{conjecture}
\begin{conjecture}\label{Conjecture:MainConjecture2} There is a semi-stable covering of $X_0(p^3)$ defined over $R_p$ such that for each supersingular elliptic curve $A$ defined over $\F_p$ there is a connected component of the supersingular locus containing (only) the following.\\
\indent (1) one component lying via $\pi_p$ over the Atkin-Lehner circle (of $A$),\\
\indent (2) two components, $E_1$ and $E_2$, which are $\pi_p^{-1}$ and $\pi_{p^2}^{-1}$ of the\\
\indent\indent too-supersingular disk (of $A$), and\\
\indent (3) $2(p+1)/i$ components lying via $\pi_p$ over the CM disks of\\
\indent\indent Conjecture \ref{Conjecture:MainConjecture1}.\\
The intersections of these components with each other and with the $6$ ordinary components of $X_0(p^3)$ are as pictured in Figure \ref{Figure:Stablep3Graph}.
\end{conjecture}

\begin{figure}\begin{center}
\includegraphics{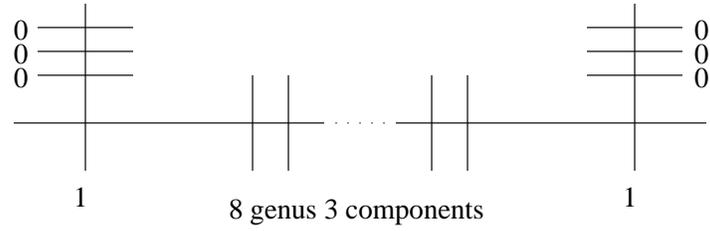}
\end{center}\caption{Conjectural Graph of $X_0(7^3)$ Semi-Stable Model, $g=26$}\label{Figure:Graph7}
\end{figure}

\begin{figure}\begin{center}
\includegraphics{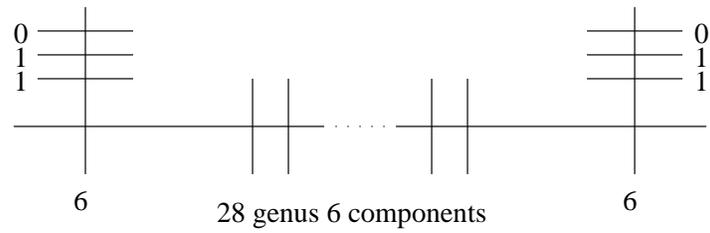}
\end{center}\caption{Conjectural Graph of $X_0(13^3)$ Semi-Stable Model, $g=184$}\label{Figure:Graph13}
\end{figure}

\begin{figure}\begin{center}
\includegraphics{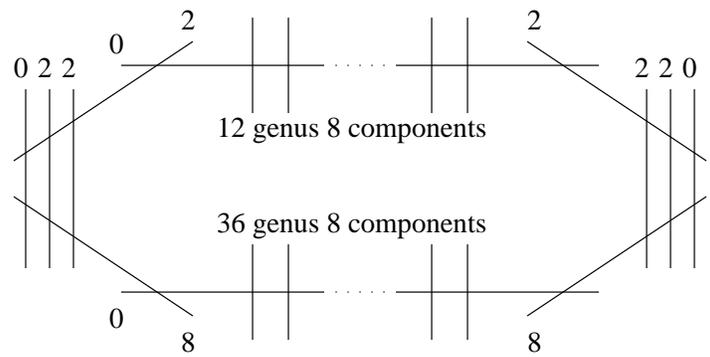}
\end{center}\caption{Conjectural Graph of $X_0(17^3)$ Semi-Stable Model, $g=417$}\label{Figure:Graph17}
\end{figure}

It is important to acknowledge that Conjecture \ref{Conjecture:MainConjecture2} is based on just one very well understood example.
Ironically, however, if we make a more precise conjecture it is easy to obtain far more corroborating data.
In particular, from looking more closely at the $X_0(125)$ example one might go so far as to guess the following.
\begin{guess} The component lying over AL is trivial, $E_1$ and $E_2$ are analytically isomorphic copies of Edixhoven's horizontal component, and each component lying over a CM disk has genus $(p-1)/2$ (and is hyperelliptic).
\end{guess}
These statements also hold for $X_0(125)$, and with the guess it becomes possible to generate complete graphs of conjectural semi-stable models for $X_0(p^3)$.
The genera of each $E_1$ and $E_2$ come directly from \cite[2.5]{E}.
The genera of the ordinary components come from \cite[2.5]{E} and the fact (from \cite{C1}) that the four nontrivial ordinary components of $X_0(p^3)$ are isomorphic copies of the two nontrivial ordinary components of $X_0(p^2)$.
Using this approach we have generated graphs for $X_0(7^3)$, $X_0(13^3)$, and $X_0(17^3)$ (Figures \ref{Figure:Graph7}, \ref{Figure:Graph13}, and \ref{Figure:Graph17}).
In each of these cases, and in fact even in the {\em general} case, it is easy to show that the total genus of the curve is at least correct.
While it is perhaps too early to make the guess an official conjecture, the corroborating data of this genus calculation seems to be very promising.

\end{document}